\documentclass[12pt]{amsart}
\usepackage{amsmath,amsfonts,amssymb,amsthm}

\pagespan{1}{?}
\usepackage{preprint}
\usepackage{graphicx}

\usepackage{psfig}

\usepackage{amsmath}

\nonstopmode
\begin{document}
\title[Experiments with moduli of quadrilaterals II]
{Experiments with moduli of quadrilaterals II}

\author[A.~Rasila]{Antti Rasila}
\email{antti.rasila@tkk.fi}
\address{Institute of Mathematics, P.O.Box 1100, FI-02015, Helsinki 
University of Technology, Finland}
\author[M.~Vuorinen]{Matti Vuorinen}
\email{vuorinen@utu.fi}
\address{Department of Mathematics, FI-20014 University of Turku, 
Finland}


\begin{abstract}
The numerical performance of the AFEM method of K.~Samuelsson is studied
in the computation of moduli of quadrilaterals.
\end{abstract}

\maketitle

\section{Introduction}

The moduli of quadrilaterals and rings are some of the fundamental
tools in geometric function theory, see \cite{Ah}, \cite{AVV},
\cite{Kuhnau:2005}, \cite{Lehto:1973}.
The purpose of this paper is to report on our experimental work
on the numerical computation of the moduli of quadrilaterals, based
on the algorithms and software of \cite{BSV} and motivated by the
geometric considerations in \cite{HVV} and \cite{DV}.
The methods considered here may be classified into two classes:
\begin{itemize}
\item[(1)] methods based on the definition of the modulus and on
conformal
mapping of the quadrilateral onto a canonical rectangle,
\item[(2)] methods based on the solution of the Dirichlet-Neumann
problem for the Laplace equation.
\end{itemize}

With the exception of a few special cases both methods lead to
extensive numerical computation. For both classes of methods
there are several options, see \cite{Gai2}, \cite{Hen}, \cite{Pap}.
Among other things, historical remarks are given in \cite{Porter}.

We study the case of a polygonal quadrilateral and the way its
modulus depends on the shape of the quadrilateral. Following the
approach of \cite{BSV} our main method is the adaptive finite
element method AFEM of Klas Samuelsson and it belongs to class (2).
We compare this method to a method of class (1), the Schwarz-Christoffel
method of L.N.~Trefethen \cite{DrTr} and its MATLAB implementation,
the SC Toolbox written by T.~Driscoll \cite{Dri}. 

In the two test cases we have used, the performance of the SC Toolbox 
was superior
to AFEM. On the other hand, the AFEM software applies also to
computation of moduli of polygonal ring domains as shown in \cite{BSV}.
AFEM also has an advantage in the problems where the quadrilateral has
a large number of vertices. This situation arises when approximating
nonpolygonal quadrilaterals (e.g. Example \ref{hyprect}).
We will  report our results also in \cite{RV}.

\section{Preliminaries}

A Jordan domain $D$ in $\C$ with marked (positively ordered) points
$z_1,z_2,z_3,z_4$ $\in \partial D$ is a quadrilateral and denoted by
$(D;z_1,z_2,z_3,z_4)\, .$ We use the canonical map of this quadrilateral
onto a rectangle $(D';1+ih,ih,0,1)$, with the vertices corresponding, to
define the modulus $h$ of a quadrilateral $(D;z_1,z_2,z_3,z_4)\,.$ The
modulus of $(D;z_2,z_3,z_4,z_1)\,$ is $1/h \, .$

We mainly study the situation
where the boundary of $D$ consists of the polygonal line segments
through $z_1,z_2,z_3,z_4$ (always positively oriented). In this case, the
modulus is denoted by $\mathrm{QM}(D;z_1,z_2,z_3,z_4)$. If the boundary
of $D$ consists of straight lines connecting the given boundary points,
we omit the domain $D$ and denote the corresponding modulus simply by
$\mathrm{QM}(z_1,z_2,z_3,z_4)$.

The following problem is known as the {\it Dirichlet-Neumann problem}.
Let $D$ be a region in the complex plane whose boundary
$\partial D$ consists of a finite number of regular Jordan
curves, so that at every point, except possibly at finitely many points
of the boundary, a normal is defined. Let $\psi$ to be a real-valued
continuous function defined on $\partial D$. Let $\partial D =A \cup B$
where $A, B$ both are unions of Jordan arcs. Find a function $u$
satisfying the following conditions:
\begin{enumerate}
\item
$u$ is continuous and differentiable in
$\overline{D}$.
\item
$u(t) = \psi(t),\qquad t \in A$.
\item
If $\partial/\partial n$ denotes differentiation in
the direction of the exterior normal, then
$$
\frac{\partial}{\partial n} u(t)=\psi(t),\qquad t\in  B.
$$
\end{enumerate}

One can express the modulus of a quadrilateral $(D; z_1, z_2, z_3, z_4)$
in terms of the solution of the Dirichlet-Neumann problem as follows.
Let $\gamma_j, j=1,2,3,4$ be the arcs of
$\partial D$ between $(z_1, z_2)\,,$ $(z_2, z_3)\,,$ $(z_3, z_4)\,,$
$(z_4, z_1),$ respectively. If $u$ is the (unique) harmonic solution of
the
Dirichlet-Neumann problem with boundary values of $u$ equal to 0 on
$\gamma_2$,
equal to 1 on $\gamma_4$ and with $\partial u/\partial n = 0$ on
$\gamma_1 \cup \gamma_3\,,$ then by \cite[p. 65/Thm 4.5]{Ah}:
\begin{equation}
\label{qmod}
\mathrm{QM}(D;z_1,z_2,z_3,z_4)=
\int_D |\triangledown
u|^2\,dm.
\end{equation}
We also have the following connection to the modulus curve family (see
e.g. \cite[pp. 158--165]{AVV}):
$\mathrm{QM}(D;z_1,z_2,z_3,z_4)=  \M(\Gamma)$,
where $\Gamma$ is the family of all curves joining $\gamma_2$ and
$\gamma_4$ in $D$.

Another approach is to use the {\it Schwarz-Christoffel formula}
to approximate the conformal mapping $f$ onto the canonical rectangle.
This formula gives an expression for a conformal map from the upper
half-plane onto the interior of a $n$-gon. Its vertices are denoted
$w_1,\ldots,w_n$, and $\alpha_1\pi,\ldots,\alpha_n\pi$ are the
corresponding interior angles. The preimages of the vertices
(prevertices) are denoted by $z_1<z_2<\ldots<z_n$. The
Schwarz-Christoffel formula for the map $f$ is
\begin{equation}
f(z) = f(z_0) + c \int_{z_0}^{z}\prod_{j=1}^{n-1}(\zeta
-z_j)^{\alpha_j-1}d\zeta,
\end{equation}
where $c$ is a (complex) constant. The main difficulty in applying this
formula is that the prevertices $z_j$ cannot, in general, be solved
for analytically. By using a M\"obius  transformation, one may choose 
three of the prevertices arbitrarily. The remaining $n-3$ prevertices 
are then obtained by solving a system of nonlinear equations. Several 
methods for solving this problem are discussed in \cite{Bis}, 
\cite{DrTr}, and \cite{DrVa}. 

The MATLAB toolbox by T. Driscoll \cite{Dri} contains a
collection of algorithms for constructing Schwarz-Christoffel maps and
computing the moduli of polygonal quadrilaterals.
The toolbox also gives an estimate for the accuracy of the numerical
approximation of the modulus.

\section{Experiments}

The solutions of the Dirichlet and the Dirichlet-Neumann problems can
be
approximated by the method of finite elements, see \cite[pp.
305--314]{Hen}, \cite{Pap}. Hence, this method can also be
used to approximate the modulus of quadrilaterals and rings.
The Dirichlet-Neumann problem can be numerically solved with AFEM
(Adaptive FEM) numerical PDE analysis package by Klas Samuelsson. This
software applies, e.g., to compute the modulus (capacity) of a bounded
ring whose boundary components are broken lines. Examples and
applications for this software are given in \cite{BSV}.

In \cite{HVV} a theoretical formula for computing $\mathrm{QM}(A,B,0,1)$ 
was given with its implementation with Mathematica. This led to a study 
of the modulus of quadrilateral in \cite{DV}. In the course of the work 
on \cite{DV}, the variation of the modulus was studied
when one of the vertices varies and others are kept fixed, and several
conjectures were formulated.
For these purposes, neither the theoretical algorithm in \cite{HVV}
nor the implemented Mathematica program based on it was no longer
adequate and we started to look for a robust program to compute
$\mathrm{QM}(A,B,0,1)\,.$
It seems that the AFEM software of Samuelsson is very efficient for this
purpose. As in \cite{romsem} we use the AFEM software of Samuelsson
for computations involving moduli of polygonal  quadrilaterals.

\begin{exmp}
Let $f(x,y)= \mathrm{QM}(x+ iy,i,0,1)- 1/\mathrm{QM}(y +ix,i, 0,1)$.
Then by \cite[p. 433]{Hen} we see that $f(x,y)\equiv 0$. Therefore
we may use this function as a measure of the accuracy of AFEM software
and SC Toolbox.
\end{exmp}

\smallskip

\centerline{
\psfig{figure=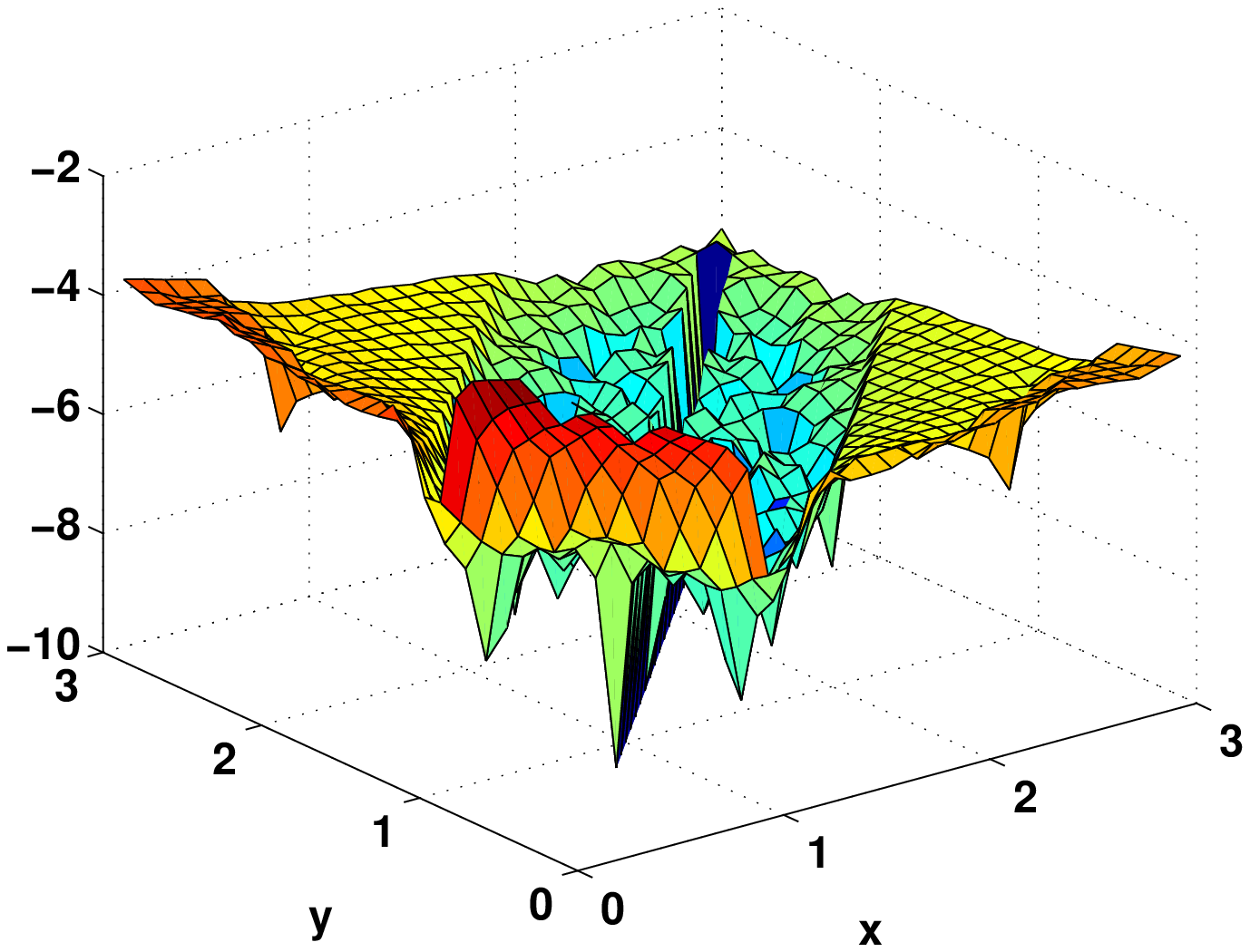,width=7cm}
\quad \psfig{figure=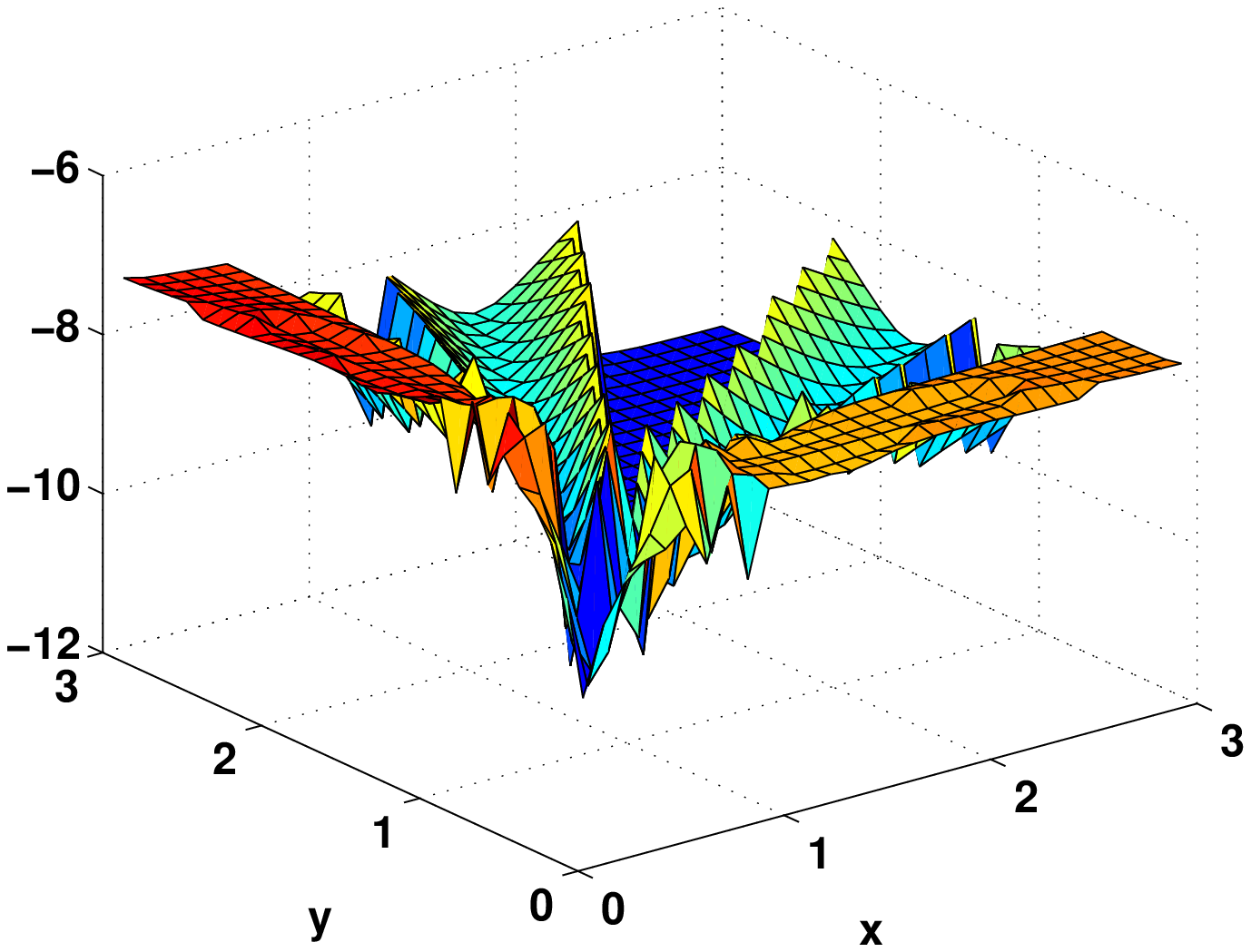,width=7.8cm}
}
\begin{capti}
Function $\log_{10}(|f(x,y)|+10^{-10})$ for $x\in(0,3]$, $y\in(0,3]$
with AFEM (left) and SC Toolbox (right).
\end{capti}

\medskip

\begin{exmp}
We study the function
$$
g(t,h) = \mathrm{QM}(1+ h e^{it},  h e^{it}, 0,1).
$$
An analytic expression for this function has been given in
\cite[2.3]{AQVV}:
\begin{equation}
g(t,h) =  \K'(r_{t/\pi})/\K(r_{t/\pi}),
\end{equation}
where
\begin{equation}
r_a = \mu_a^{-1}\bigg(\frac{\pi h}{2\sin(\pi a)}\bigg),
\text{ for }0<a\leq 1/2,
\end{equation}
and the decreasing homeomorphism $\mu_a\colon (0,1)\to(0,\infty)$
is defined by
\begin{equation}
\mu_a(r)\equiv \frac{\pi}{2\sin(\pi
a)}\,\frac{F(a,1-a;1;1-r^2)}{F(a,1-a; 1;
r^2)}.
\end{equation}
Here
$$
F(a,b;c;z) = {}_2F_1(a,b;c;z) \equiv \sum_{n=0}^\infty
\frac{(a,n)(b,n)}{(c,n)}\, \frac{z^n}{n!}, \qquad |z|<1,
$$
is the {\it Gaussian hypergeometric function},
$$
(a,n)\equiv a(a+1)(a+2)\ldots(a+n-1),\qquad
(a,0)=1\text{ for }a\neq 0,
$$
is the {\it shifted factorial function}, and the
{\it elliptic integrals} $\K(r),\K'(r)$
are defined by
$$
\K(r)=\frac{\pi}{2} F(1/2,1/2;1; r^2),
\qquad
\K'(r)=\K(r'),\text{ and }r'= \sqrt{1-r^2}.
$$


\centerline{
\psfig{figure=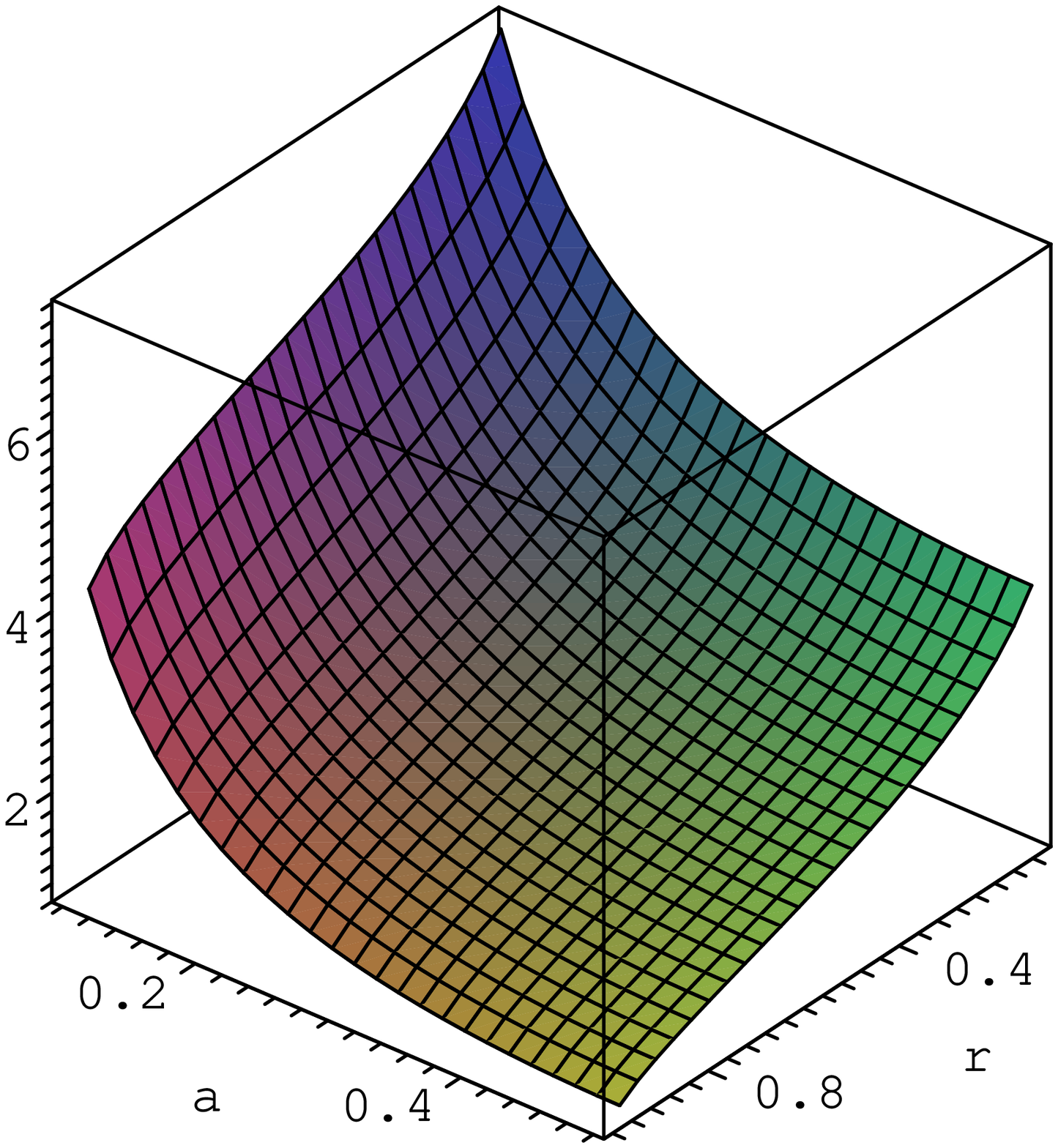,width=8cm}
}
\vskip -.3cm
\begin{center}
\begin{capti}
Function $\mu_a(r)$.\vskip .5cm
\end{capti}
\end{center}

The function $g(t,h)$ is the modulus of the parallelogram with opposite
sides  $1$ and $h$, respectively, and we see that there are three cases
$h \in (0,1)$, $h=1$ and $h>1$. In the first case the function is
monotone increasing with respect to $t \in (0, \pi/2)$, in the second case the
function $g(t,1)\equiv 1$ is constant and in the third case decreasing.

\medskip

\centerline{
\psfig{figure=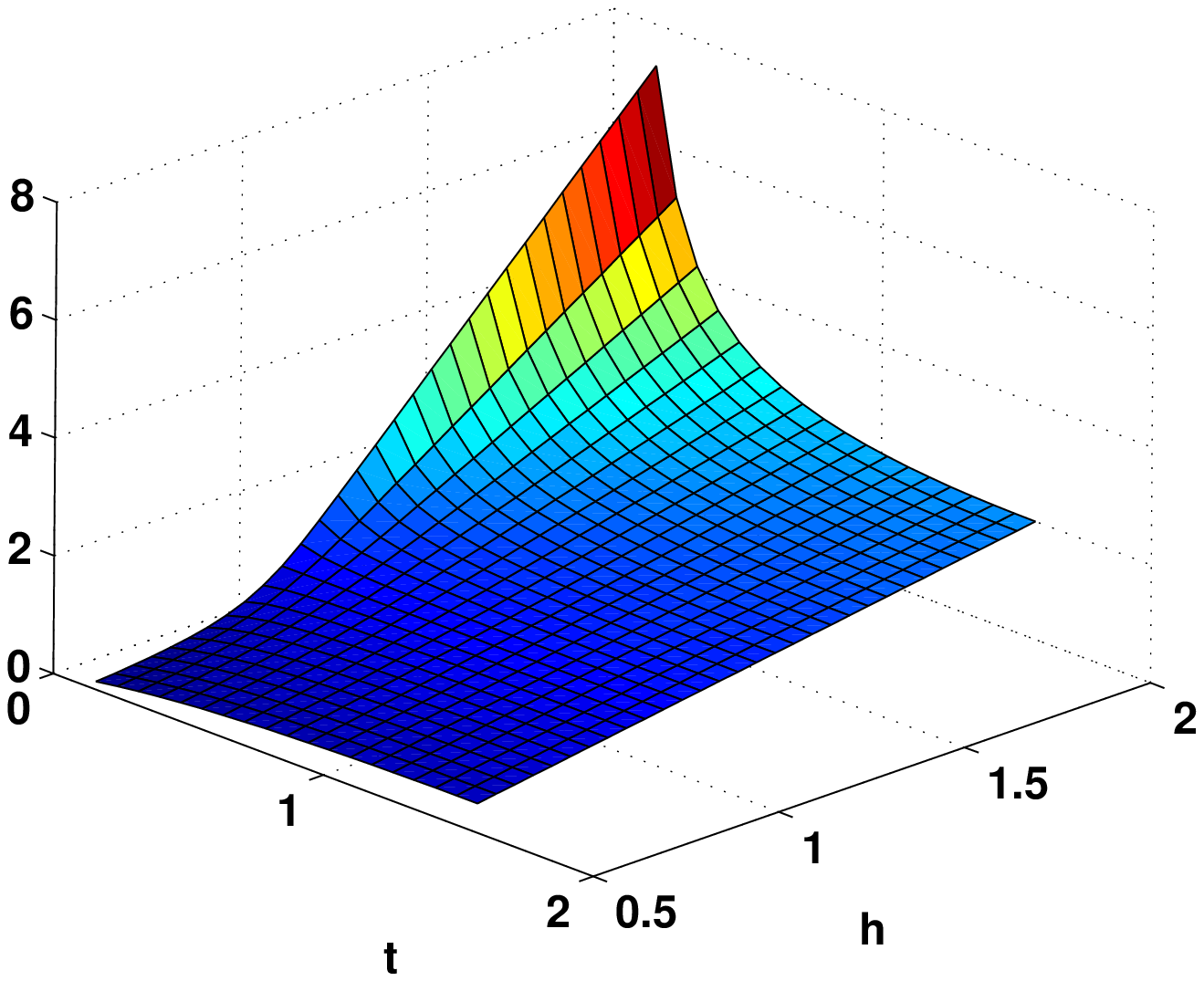,width=7cm}
\psfig{figure=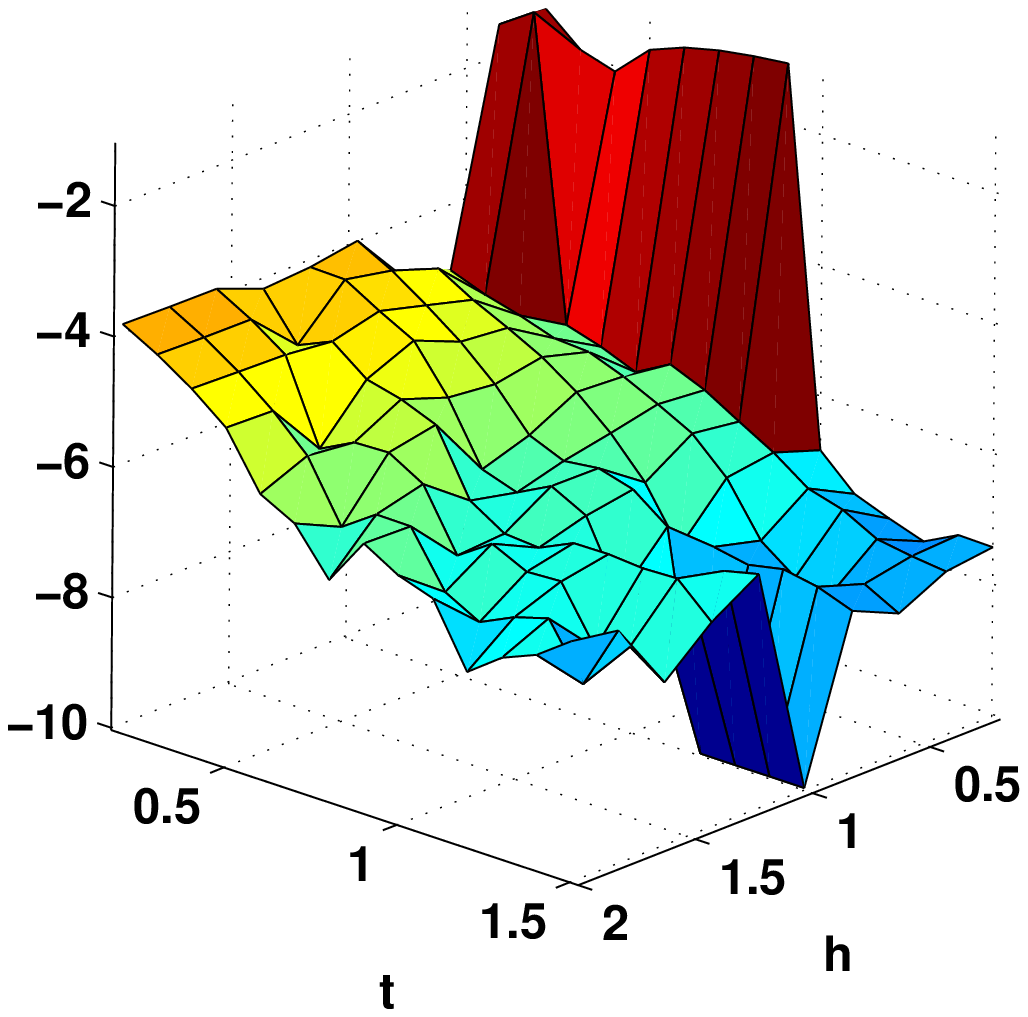,width=8cm}
}
\begin{capti}
Function $g(t,h)$ for $t\in(0,\pi/2)$ and $h\in[1/2,2]$ (left), and
the error estimate
$\log_{10}(|g_{\mathrm{exact}}(t,h)-g_{\mathrm{numer}}(t,h)|+10^{-10})$
for
the function $g(t,h)$ (right).
\end{capti}

\medskip

\centerline{
\psfig{figure=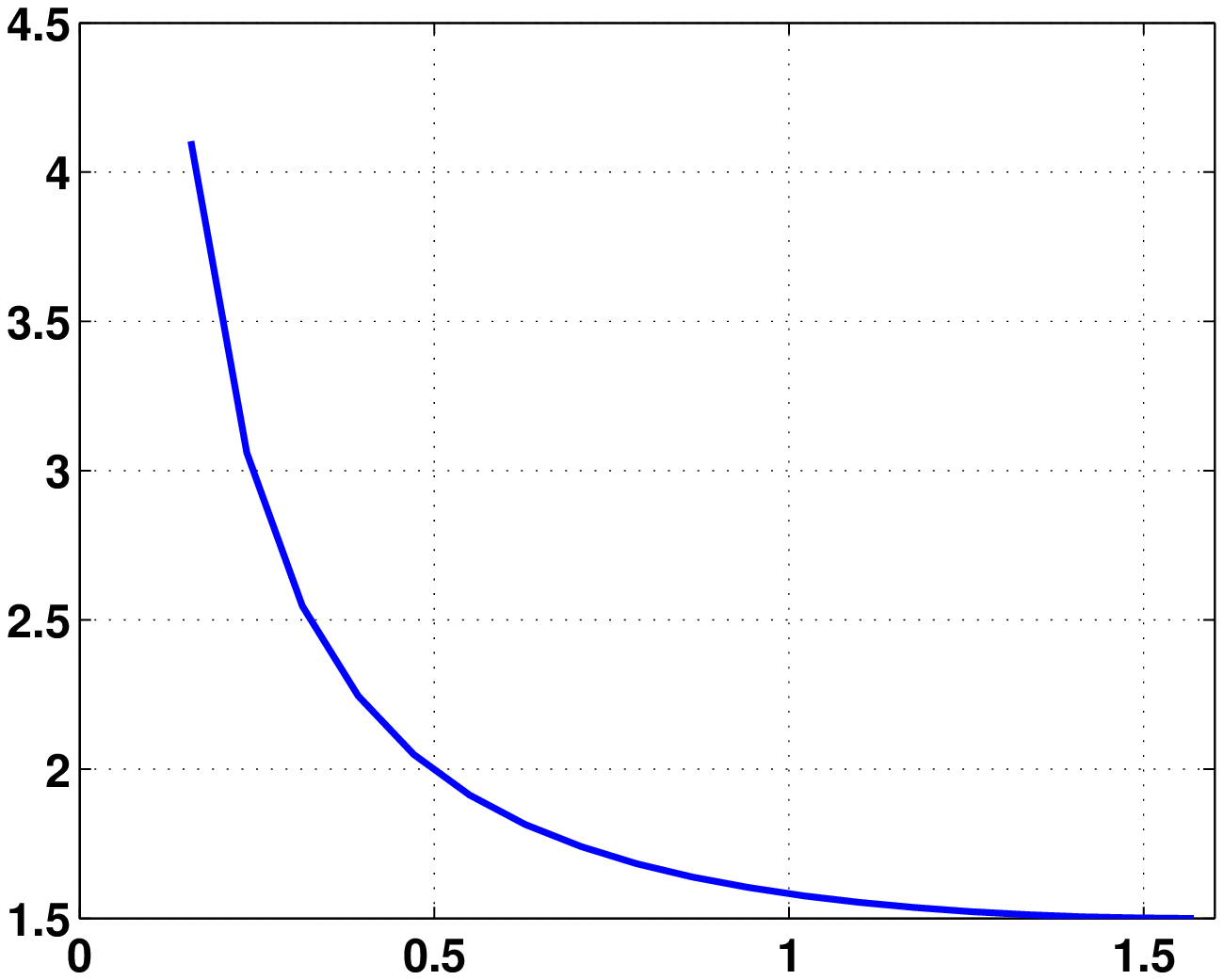,width=10.5cm}
}
\begin{center}
\vskip-0.5cm
\begin{capti}
Function $g(t,1.5)$ for $t\in(0,\pi/2)$.
\end{capti}
\end{center}
\end{exmp}


\medskip

\begin{exmp}
\label{mrem}
The modulus $\mathrm{QM}(1+i|a-1|,i|b|,0,1)$ has an
analytic expression if $|a-1|= h= |b|+1\,.$
Bowman \cite[pp.\ 103-104]{Bo} (see also \cite[p.197]{Bur}) gives a
formula for the conformal modulus of the quadrilateral with vertices
$1+hi$,  $(h-1)i$, $0$, and $1$ when $h>1$ as $M(h)\equiv\K(r)/\K(r')$
where
$$
r = \bigg(\frac{t_1-t_2}{t_1+t_2}\bigg)^2\, , \quad
t_1 = \mu_{1/2}^{-1}\left(\frac{\pi}{2c}\right) \, , \quad
t_2 = \mu_{1/2}^{-1}\left(\frac{\pi c}{2}\right)\, , \quad
c = 2h - 1 \, .
$$
Therefore, the quadrilateral can be conformally mapped onto the
rectangle  $1+iM(h)$, $i M(h)$, $0$, $1$, with the vertices
corresponding
to each other. It is clear that $h-1 \le M(h) \le h \,.$
The formula
$$M(h)= h + c +O(e^{-\pi h}),\qquad  c=-1/2- \log
2/\pi\approx -0.720636\,,
$$
is given in \cite{PS}. As far as we know there is neither
an explicit nor asymptotic formula
for the case when the angle $\pi/4$ of the trapezoid is equal
to $\alpha \in (0, \pi/2) \,.$
We compute the modulus $\mathrm{QM}(ih,i(h-1),0,1)$ by using
the square frame capacity formula, AFEM and Schwarz-Christoffel Toolbox.

\medskip
{\rm \Small
$$
\begin{array}{|c|l|l|l|l|l|l|}
\hline h   &  \textbf{AFEM }  & \textbf{SC} & \textbf{Accuracy} &
\textbf{Bowman} & \textbf{Error} & \textbf{Error/SC}\\
&&&\textbf{/SC}&&\textbf{/AFEM}&\\
\hline 1.1  &  0.3403159  &  0.3403135 & 1.787\cdot 10^{-8}  &
0.3403135 &
2.41655\cdot 10^{-6}  &  1.57002\cdot 10^{-9}\\
1.2  &  0.4614938  &  0.4614926 & 1.734\cdot 10^{-8}  &  0.4614926
&
1.20727\cdot 10^{-6}    &2.98441\cdot 10^{-9}\\
1.3  &  0.5704380  &  0.5704374 & 5.310\cdot 10^{-8}  &  0.5704374
&
5.83493\cdot 10^{-7}    &4.59896\cdot 10^{-9}\\
1.4  &  0.6747519  &  0.6747518 & 1.046\cdot 10^{-7}  &  0.6747518
&
8.83554\cdot 10^{-8}    &6.24458\cdot 10^{-9}\\
1.5  &  0.7769433  &  0.7769434 & 2.408\cdot 10^{-8}  &  0.7769434
&
1.10607\cdot 10^{-7}    &3.39673\cdot 10^{-9}\\
1.6  &  0.8780836  &  0.8780838 & 1.920\cdot 10^{-9}  &  0.8780838
&
1.53305\cdot 10^{-7}    &8.10543\cdot 10^{-10}\\
1.7  &  0.9786840  &  0.9786842 & 5.439\cdot 10^{-10}  &  0.9786842
&
2.41392\cdot 10^{-7}    &2.02109\cdot 10^{-10}\\
1.8  &  1.0790020  &  1.0790024 & 2.102\cdot 10^{-10}  &  1.0790024
&
4.03325\cdot 10^{-7}    &4.94438\cdot 10^{-11}\\
1.9  &  1.1791710  &  1.1791715 & 6.225\cdot 10^{-11}  &  1.1791715
&
5.22481\cdot 10^{-7}    &1.20739\cdot 10^{-11}\\
2.0  &  1.2792610  &  1.2792616 & 1.536\cdot 10^{-11}  &  1.2792616
&
5.71171\cdot 10^{-7}    & 2.97451\cdot 10^{-12}\\
\hline
\end{array}
$$
}

\begin{tabti}
Error estimates for AFEM and SC Toolbox with the square frame capacity
formula. The accuracy of the estimate given by SC Toolbox is also 
consistent with the experiment.
\end{tabti}

\end{exmp}

\begin{exmp}
\label{hyprect}
Let $Q$ be the quadrilateral whose sides are defined by two circular
arcs in the upper and lower half plane, perpendicular to the unit
circle at the points $e^{i\theta}$, $e^{i(\pi-\theta)}$,
$e^{i(\theta-\pi)}$, $e^{-i\theta}$ as well
as by the two circular arcs through $r$,  $i$, $-i$ and $-r$, $i$, $-i$,
see Figure 6. If $a,b,c,d$ are the points of intersection of
these four circular arcs in the II, III, IV and I quadrants
respectively, then $Q=(Q; a,b,c,d)$ defines a quadrilateral in the unit
disk with
$$\mathrm{QM}(Q; a,b,c,d)= (\pi- 2 \beta)/\rho,
$$
where
$$
\rho=2\log\frac{1+u}{1-u}, \quad
\beta=\mathrm{arccot}\frac{2r}{1-r^2},\text{ and }u=\tan(\theta/2).
$$

\end{exmp}


\centerline{
\psfig{figure=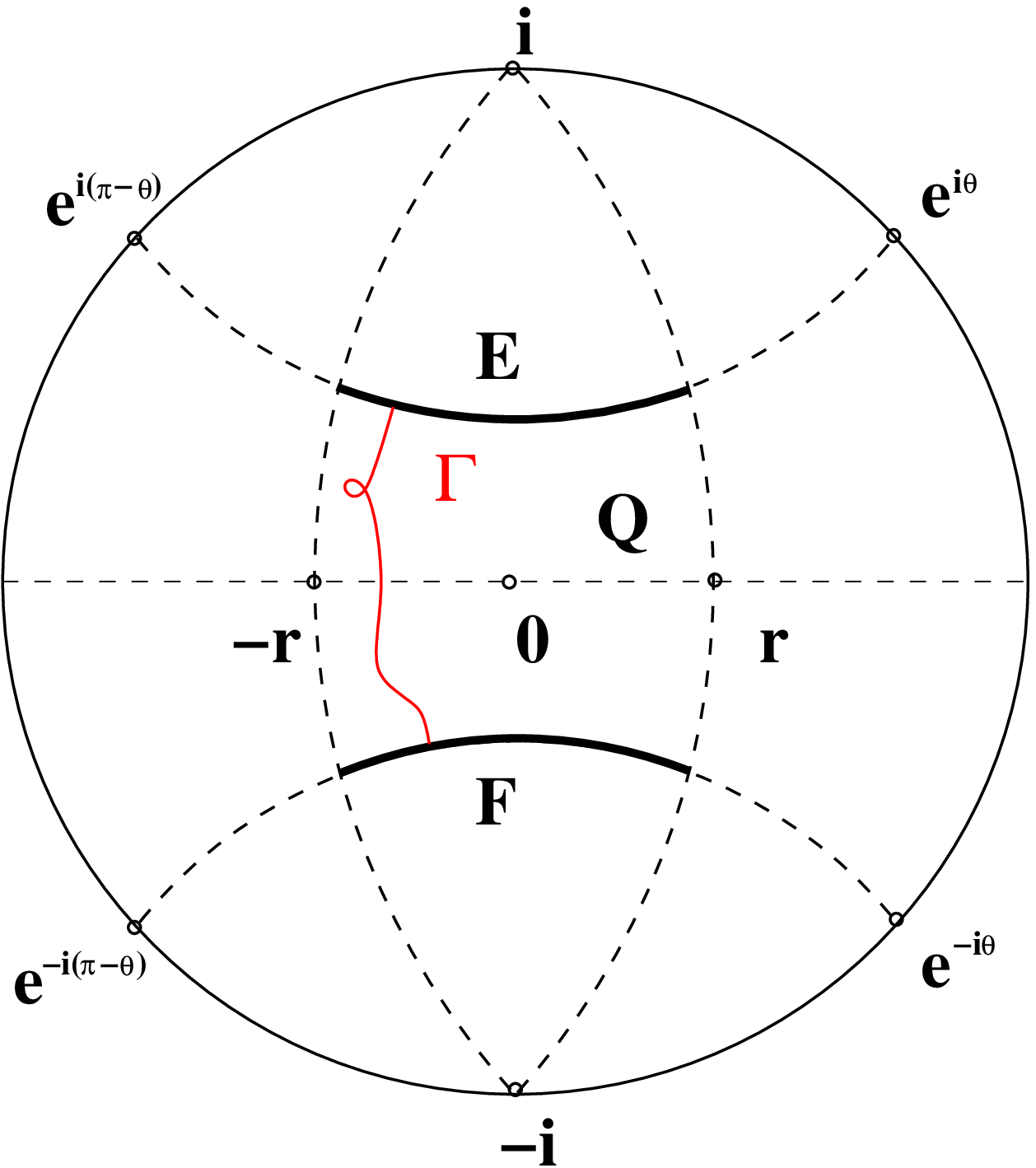,width=5.7cm}
}
\begin{center}
\begin{capti}
\label{below}
The circular quadrilateral $Q$.
\end{capti}
\end{center}

\medskip

$$
\begin{array}{|c|l|l|l|}
\hline
     \theta & \textbf{AFEM} &  \textbf{Exact}&\textbf{Error}\\
\hline
0.10  &  7.592357  &  7.597433 &  5.076\cdot 10^{-3}\\
0.15  &  5.056044  &  5.054357 &  1.687\cdot 10^{-3}\\
0.20  &  3.784480  &  3.779611 &  4.869\cdot 10^{-3}\\
0.25  &  3.010221  &  3.012175 &  1.954\cdot 10^{-3}\\
0.30  &  2.497983  &  2.498368 &  3.849\cdot 10^{-4}\\
0.35  &  2.130426  &  2.129465 &  9.616\cdot 10^{-4}\\
\ldots & \ldots & \ldots & \ldots \\
1.00  &  0.620785  &  0.620631 &  1.531\cdot 10^{-4}\\
1.05  &  0.575335  &  0.575402 &  6.645\cdot 10^{-5}\\
1.10  &  0.533529  &  0.533010 &  5.185\cdot 10^{-4}\\
1.15  &  0.493114  &  0.492934 &  1.805\cdot 10^{-4}\\
1.20  &  0.454678  &  0.454689 &  1.114\cdot 10^{-5}\\
\hline
\end{array}
$$
\begin{center}
\begin{tabti}
The modulus of the quadrilateral $Q$ for $r=0.4$.
\end{tabti}
\end{center}

\end{document}